\documentclass[12pt, reqno]{amsart}

\usepackage{amssymb, amsthm}
\usepackage{fullpage}

\numberwithin{equation}{section}

\newtheorem{lemma}{Lemma}[section]

\newtheorem{thm}[lemma]{Theorem}

\theoremstyle{remark}
\newtheorem{remark}{Remark}

\newcommand{\imi}{\mathrm i}
\newcommand{\eps}{\varepsilon}
\newcommand{\modulo}[1]{\;\mathrm{mod}\;#1}

\begin{document}

\title[Mean Values of Dirichlet Polynomials]{Mean Values of Dirichlet Polynomials and Applications to Linear Equations with Prime Variables}
\author[S.K.K. Choi]{Stephen Kwok-Kwong Choi}
\address{Department of Mathematics and Statistics \\
Simon Fraser University \\
Burnaby, British Columbia \\
Canada V5A 1S6}
\email{kkchoi@math.sfu.ca}
\thanks{Research of Stephen Choi was supported by NSERC of Canada.}
\author[A.V. Kumchev]{Angel V. Kumchev}
\address{Department of Mathematics \\
1 University Station C1200 \\
The University of Texas at Austin \\
Austin, TX 78712 \\
U.S.A.}
\email{kumchev@math.utexas.edu}
\date{Draft from \today}

\maketitle

\section{Introduction}

In this paper we study Dirichlet polynomials of the form
\begin{equation}\label{i.1}
  D(s, \chi) = \sum_{n \le N} a_n \chi(n) n^{-s}
\end{equation}
where $\chi(n)$ is a Dirichlet character, $s = \sigma + \imi t$ is a complex variable, and $a_n$ are (complex) coefficients. Such Dirichlet polynomials are an important tool in multiplicative number theory and there is a vast literature on the subject. In particular, one often needs estimates for mean values of the form
\[
  \sum_{\chi \in \mathcal H} \int_{-T}^T \bigg| \sum_{N < n \le 2N} \Lambda(n)\chi(n)n^{-\imi t} \bigg| \, dt,
\]
where $\Lambda(n)$ is the von Mangoldt function and the outer summation is over some family of characters, possibly to various moduli. Our main result is Theorem \ref{th1} below, which deals with the most common types of such averages.

Let $m \ge 1$, $r \ge 1$, and $Q \ge r$. We consider a set $\mathcal H(m, r, Q)$ of characters $\chi = \xi\psi$ modulo~$mq$, where $\xi$ is a character modulo $m$ and $\psi$ is a primitive character modulo~$q$, with $r \le q \le Q$, $r \mid q$, and $(q, m) = 1$. Our result is as follows.

\begin{thm}\label{th1}
  Let $m \ge 1$, $r \ge 1$, $Q \ge r$, $T \ge 2$, $N \ge 2$, and $\mathcal H(m, r, Q)$ be a set of characters as described above. Then 
  \begin{equation}\label{i.2}
    \sum_{\chi \in \mathcal H(m, r, Q)} \int_{-T}^T 
    \bigg| \sum_{N < n \le 2N} \Lambda(n)\chi(n)n^{-\imi t} \bigg| \, dt
    \ll \big( N + HN^{11/20} \big) L^C,
  \end{equation}
  where $C$ is an absolute constant,
  \[
    H = mr^{-1}Q^2T \qquad \text{and} \qquad L = \log HN.
  \]
\end{thm} 

\begin{remark}
  A possible choice for $C$ is $C = 1100$, and we have organized the proof as to make this obvious. On the other hand, we have spent no effort to optimize our estimates in that regard, because it is clear that our method will never yield a result with a ``respectable'' value of $C$, such as $C = 10$, or even $C = 100$.
\end{remark}

\begin{remark}
Under the Generalized Riemann Hypothesis (GRH), we have
\[
  \sum_{\chi \in \mathcal H(m, r, Q)} \int_{-T}^T 
  \bigg| \sum_{N < n \le 2N} \Lambda(n)\chi(n)n^{-\imi t} \bigg| \, dt
  \ll NL + HN^{1/2}L^2,
\]
where the term $NL$ on the right side occurs only when the set $\mathcal H(m, r, Q)$ contains a principal character. In contrast, because Theorem \ref{th1} is derived from a general result on bilinear forms (see Theorem \ref{th2.1} below), the first term on the right side of \eqref{i.2} occurs independent of the presence of a principal character in $\mathcal H(m, r, Q)$.
\end{remark}

Using Theorem \ref{th1}, we can make progress in an additive problem with prime variables. Consider the linear diophantine equation 
\begin{equation}\label{2.1}
  a_1p_1 + a_2p_2 + a_3p_3 = b
\end{equation}
where $a_1, a_2, a_3, b$ are integers with $a_1a_2a_3 \ne 0$ and $p_1, p_2, p_3$ are prime unknowns. Our goal is to prove the existence of solutions of \eqref{2.1} which do not grow too rapidly as $B = \max\{ |a_1|, |a_2|, |a_3| \} \to \infty$. This problem was first raised and investigated by Baker \cite{Ba} and was later settled, at least qualitatively, by M.C. Liu and Tsang~\cite{LT}. A necessary condition for the solubility of \eqref{2.1} is
\begin{equation}\label{C1}
  a_1 + a_2 + a_3 \equiv b \pmod 2.
\end{equation}
Without loss of generality, we may assume that
\begin{equation}\label{C2}
  (a_1, a_2, a_3) = (b, a_i, a_j) = 1, \qquad 1 \le i < j \le 3.
\end{equation}
Liu and Tsang \cite{LT} proved the following result.

\begin{thm}\label{thm LT}
  Suppose that $a_1,a_2,a_3,b$ are integers such that $a_1a_2a_3 \ne 0$ and conditions \eqref{C1} and \eqref{C2} hold. Then there exists an absolute constant $A > 0$ such that
  \begin{enumerate}
    \item[(i)]  if $a_1,a_2,a_3$ are all positive, then \eqref{2.1} has solutions in primes whenever $b \gg B^A$;
    \item[(ii)] if $a_1,a_2,a_3$ are not all of the same sign, then \eqref{2.1} has solutions in primes satisfying
    \begin{equation}\label{2.4}
      |a_j|p_j \ll |b| + B^A.
    \end{equation}
  \end{enumerate}
\end{thm}

It is not difficult to see that one cannot take the exponent $A$ above arbitrarily small, so it remains to estimate the best possible value of $A$. The first numerical upper bound for $A$ was obtained by Choi \cite{C}, who showed that $A \le 4190$. This bound was subsequently reduced to $A \le 45$ by M.C. Liu and Wang \cite{LW} and to $A \le 38$ by Li~\cite{Li}. Furthermore, Choi, M.C. Liu, and Tsang~\cite{CLT} showed that under GRH one has $A \le 5 + \eps$ for any fixed $\eps > 0$.

Recently, J.Y. Liu and Tsang \cite{jLT} showed that when condition \eqref{C2} is replaced by the somewhat more restrictive
\begin{equation}\label{2.5}
  (a_1, a_2) = (b, a_i)=1, \qquad 1 \le i < j \le 3,
\end{equation}
then one can take (essentially) $A = 17/2$. In the last section of this paper, we obtain the following improvement on their result, thus reducing the value of $A$ further to $A = 20/3$.

\begin{thm}\label{th3}
  Suppose that $a_1,a_2,a_3,b$ are integers such that $a_1a_2a_3 \ne 0$ and conditions \eqref{C1} and \eqref{2.5} hold.
  \begin{enumerate}
    \item[(i)] If $a_1,a_2,a_3$ are all positive, then \eqref{2.1} has solutions in primes whenever 
    \[
      b \gg (a_1a_2a_3)^{20/9}B(\log B)^{26}.
    \]
    \item[(ii)] If $a_1,a_2,a_3$ are not all of the same sign, then \eqref{2.1} has solutions in primes satisfying
    \[
      |a_j|p_j \ll |b| + (a_1a_2a_3)^{20/9}B(\log B)^{26}.
    \]
  \end{enumerate}
\end{thm}

\begin{remark}
The proof of Theorem \ref{thm LT} uses the circle method and the Deuring--Heilbronn phenomenon to treat the major arcs, which need to be taken significantly larger than in classical applications. Under the condition \eqref{C2} in place of \eqref{2.5}, one can show that the possible existence of Siegel zeros does not have special influence, and hence the Deuring--Heilbronn phenomenon can be avoided (see \cite[Lemma 3.1]{jLT}). As a result, better results can be obtained without recourse to the heavy numerical computations needed in \cite{C, Li, LW}.
\end{remark}

\section{Mean values of products of Dirichlet polynomials}

We derive Theorem \ref{th1} from mean-value estimates for products of Dirichlet polynomials of the form
\begin{equation}\label{1.1}
  F(s, \chi) = \prod_{i = 1}^3 \bigg\{ \sum_{N_i < n \le N_i'} b_i(n) \chi(n) n^{-s} \bigg\}.
\end{equation}
We assume that $1 \le N_i < N_i' \le 2N_i$ and $X = N_1N_2N_3 \ge 10$. We also assume that the coefficients $b_j(n)$ are subject to
\begin{equation}\label{1.2}
  |b_1(n)| \le \tau_{\kappa}(n), \qquad |b_2(n)| \le \tau_{\nu}(n), \qquad |b_3(n)| \le 1
\end{equation}
for some integers $\kappa, \nu \ge 2$. Here, $\tau_{\kappa}(n)$ denotes the $\kappa$-fold divisor function. The main result of this section is the following theorem.

\begin{thm}\label{th2.1}
  Suppose that $\mathcal H(m, r, Q)$ is a set of characters as in Theorem \ref{th1} and $F(s, \chi)$ is a Dirichlet polynomial as above. Also, suppose that either
  \begin{enumerate}
    \item  [(i)] $\max(N_1, N_2) \ll X^{11/20}$ and $b_3(n) = 1$ for all $n \le 2N_3$, or
    \item [(ii)] $\max(N_1, N_2) \ll X^{11/20}$ and $N_3 \ll X^{8/35}$.
  \end{enumerate}
  Then
  \begin{equation}\label{1.3}
    \sum_{\chi \in \mathcal H(m, r, Q)} \int_{-T}^T \big| F(\imi t, \chi) \big| \, dt
    \ll \big( X + HX^{11/20} \big) L^{c(\kappa, \nu)},
  \end{equation}
  where $c(\kappa, \nu) = 3\max(\kappa^2, \nu^2) + \kappa + \nu + 20$, $H = mr^{-1}Q^2T$, and $L = \log 2HX$.
\end{thm}

The main tool in the proof of Theorem \ref{th2.1} are bounds for the cardinality of a well-spaced set of points at which a Dirichlet polynomial of the form \eqref{i.1} is large. In this context, a ``point'' is an ordered pair $(t, \chi)$, where $t$ is a real number such that $|t| \le T$ and $\chi$ is a character from $\mathcal H(m, r, Q)$. We say that the points $(t_1, \chi_1), \dots, (t_R, \chi_R)$ are \emph{well-spaced} if $|t_i - t_j| \ge 1$ whenever $\chi_i = \chi_j$ and $i \ne j$. 

\begin{lemma}\label{l2.2}
  Suppose that $(t_1, \chi_1), \dots, (t_R, \chi_R)$ are well-spaced and that for all $j = 1, \dots, R$,
  \[
    \bigg| \sum_{n \le N} a_n\chi_j(n)n^{-\imi t_j} \bigg| \ge V.
  \]
  Then
  \[
    R \ll \big( NV^{-2} + H \min \big\{ V^{-2}, NG^2V^{-6} \big\} \big) GL^{18},
  \]
  where $L = \log 2HN$ and
  \[
    G = \sum_{n \le N} |a_n|^2.
  \]
\end{lemma}
 
\begin{proof}
  When $r = 1$, the lemma is a direct consequence of \cite[Theorem 9.16]{IK} and \cite[Theorem 9.18]{IK}. When $r > 1$, we need respective modifications of those results. The modifications, however, are straightforward because of the following observations:
  \begin{itemize}
    \item the trivial bound for the cardinality of $\mathcal H(m, r, Q)$ is $|\mathcal H(m, r, Q)| \ll mr^{-1}Q^2$;
    \item if $1 \le q_1, q_2 \le Q$ and $r \mid (q_1, q_2)$, then $[q_1, q_2] \le r^{-1}Q^2$. 
  \end{itemize}
\end{proof}

\begin{lemma}\label{l2.3}
  Let $N < M \le cN$ and define
  \begin{equation}\label{Ddef}
    D(s, \chi) = \sum_{N < n \le M} \chi(n)n^{-s}.
  \end{equation}
  Suppose that $(t_1, \chi_1), \dots, (t_R, \chi_R)$ are well-spaced and that $|t_j| \ge N$ whenever $\chi_j$ is principal. Then
  \begin{equation}\label{4PM}
    \sum_{j = 1}^R \left| D(\imi t_j, \chi_j) \right|^4 \ll HN^2L^{10}.
  \end{equation}
\end{lemma}

\begin{proof}
  Without loss of generality we may assume that the distances from $M$ and $N$ to $\mathbb Z$ equal $1/2$. For any character $\chi \in \mathcal H(m, r, Q)$, Perron's formula (see \cite[Proposition 5.54]{IK}) yields
  \[
    D(\imi t, \chi) = \frac 1{2\pi \imi } \int_{\alpha - \imi T_1}^{\alpha + \imi T_1} L(\imi t + w, \chi) \frac {M^w - N^w}{w} \, dw + O(1),
  \]
  where $T_1 = 10HN$ and $\alpha = 1 + (\log T_1)^{-1}$. The integrand is holomorphic everywhere except possibly at $w = 1 - \imi t$, where $L(\imi t + w, \chi)$ has a simple pole if $\chi$ is principal. Thus, we can move the integration to the contour $\mathfrak C$ consisting of the other three sides of the rectangle with vertices $1/2 \pm \imi T_1, \alpha \pm \imi T_1$. By the convexity bound
  \[
    L(\sigma + \imi t, \chi) \ll (mq(|t| + 2))^{(1 - \sigma)/2 + \eps}
    \qquad (0 \le \sigma \le 1),
  \] 
  the integrals over the horizontal parts of $\mathfrak C$ contribute at most
  \[
    \sup_{1/2 \le \sigma \le \alpha}\big\{ T_1^{-1}N^{\sigma} (mqT_1)^{(1 - \sigma)/2 + \eps} \big\} \ll 1.
  \]
  Also, the residue at $w = 1 - \imi t$ is $\ll \delta_{\chi}NL(1 + |t|)^{-1}$, where $\delta_{\chi}$ is $1$ or $0$ according as $\chi$ is principal or not. Hence, for any point $(t_j, \chi_j)$, $j = 1, \dots, R$, we have
  \begin{align*}
    D(\imi t_j, \chi_j) &\ll N^{1/2}\int_{-T_1}^{T_1} \big| L( 1/2 + \imi (t_j + u), \chi_j) \big| \frac {du}{1 + |u|} + \frac {\delta_{\chi_j}NL}{1 + |t_j|} + 1\\
    &\ll N^{1/2}\int_{-T_1}^{T_1} \big| L( 1/2 + \imi (t_j + u), \chi_j) \big| \frac {du}{1 + |u|} + L,
  \end{align*}
  where the last inequality uses the hypothesis on points $(t_j, \chi_j)$ with principal characters. Appealing to H\"older's inequality, we derive the estimate
  \begin{align*}
    |D(\imi t_j, \chi_j)|^4 &\ll N^2L^3
    \int_{-T_1}^{T_1} \big| L(1/2 + \imi (t_j + u), \chi_j) \big|^4 \frac {du}{1 + |u|} + L^4 \\
    &\ll N^2L^3
    \int_{-2T_1}^{2T_1} \big| L(1/2 + \imi u, \chi_j) \big|^4 \frac {du}{1 + |u - t_j|} + L^4,
  \end{align*}
  whence
  \[
    \sum_{j = 1}^R |D(\imi t_j, \chi_j)|^4  \ll N^2L^3 \sum_{j = 1}^R \int_{-2T_1}^{2T_1} 
    \big| L(1/2 + \imi u, \chi_j) \big|^4 \frac {du}{1 + |u - t_j|} + RL^4.
  \]
  This suffices, because
  \begin{align*}
    & \sum_{j = 1}^R \int_{-2T_1}^{2T_1} \big| L(1/2 + \imi u, \chi_j) \big|^4 \frac {du}{1 + |u - t_j|} \\
    \ll \; & \sum_{\chi \in \mathcal H(m, r, Q)} \int_{-2T_1}^{2T_1} \big| L(1/2 + \imi u, \chi) \big|^4 \bigg\{ \sum_{ \substack{ j = 1\\ \chi_j = \chi}}^R \frac 1{1 + |u - t_j|} \bigg\} du \\
    \ll \; & TL \sum_{\chi \in \mathcal H(m, r, Q)} \int_{-2T_1}^{2T_1} \big| L(1/2 + \imi u, \chi) \big|^4 \frac {du}{T + |u|} \ll HL^7,
  \end{align*}
  where the final step uses the estimate for the fourth power moment of $L(s, \chi)$ (see \cite[Theorem 10.1]{Mo}).
\end{proof}

\begin{proof}[Proof of Theorem \ref{th2.1}]
  Define the Dirichlet polynomials
  \[
    F_i(s, \chi) = \sum_{N_i < n \le N_i'} b_i(n) \chi(n) n^{-s} \qquad (i = 1, 2, 3).
  \]  
  The proof is divided into several steps.
  
  \medskip
  
  \paragraph{\emph{Step 1}}
  First, we dispense with some technical difficulties caused by the principal character $\chi_0$ modulo $m$ (if present in $\mathcal H(m, r, Q)$) when we argue under hypothesis (i). By the properties of the M\"obius function,  
  \[
    |F_3(\imi t, \chi_0)| \le \sum_{d \mid m} \bigg| \sum_{N_3 < dn \le N_3'} n^{-\imi t} \bigg|.
  \]
  Hence,
  \begin{equation}\label{2.1.a}
    \int_{-T}^T \big| F(\imi t, \chi_0) \big| dt 
    \ll L\sum_{d \mid m} \int_{-T}^T \big| G_d(\imi t) \big| \, dt,
  \end{equation}
  where 
  \[
    G_d(s) = \sum_{N_1 < n_1 \le N_1'} \sum_{N_2 < n_2 \le N_2'} \sum_{M_d < n_3 \le M_d'} 
    \tilde b_1(n_1) \tilde b_2(n_2) (n_1n_2n_3)^{-s},
  \]
  with $M_d < M_d' \le 2N_3/d$ and coefficients subject to
  \[
    |\tilde b_1(n)| \le |b_1(n)|, \qquad |\tilde b_2(n)| \le |b_2(n)|.
  \]
  We now recall the well-known estimates (see \cite[(1.80)]{IK} and \cite[Corollary 8.11]{IK})
  \[
    \sum_{n \le x} \tau_{\kappa}(n)^{\nu} \ll x(\log x)^{\kappa^{\nu} - 1}
  \]
  and
  \[
    \sum_{N < n \le 2N} n^{-it} \ll N(1 + |t|)^{-1} \qquad (|t| < N).
  \]
  Using the former bound to estimate the sums over $n_1$ and $n_2$ and the latter to estimate the sum over $n_3$, we get
  \begin{align}\label{2.1.b}
    \sum_{d \mid m} \int_{-M_d}^{M_d} \big| G_d(\imi t) \big| \, dt
    &\ll N_1N_2L^{\kappa + \nu -2} \sum_{d \le 2N_3} \int_{-M_d}^{M_d} \frac {M_d' \, dt}{1 + |t|} \notag\\
    &\ll N_1N_2L^{\kappa + \nu - 1} \sum_{d \le 2N_3} N_3d^{-1} \ll XL^{\kappa + \nu}. 
  \end{align}
  On the other hand, for each $d \mid m$ such that $M_d < T$, the estimates in Steps 4 and~5 below yield
  \begin{equation}\label{2.1.c}
    \int_{M_d \le |t| \le T} |G_d(\imi t)| \, dt \ll \big( X_d + TX_d^{11/20} \big)L^{c_0},
  \end{equation}
  where $c_0 = c_0(\kappa, \nu) = 3\max(\kappa^2, \nu^2) + \kappa + \nu + 15$ and $X_d = Xd^{-1}$. Thus,
  \begin{align}\label{2.1.d}
    \sum_{d \mid m} \int_{M_d \le |t| \le T} |G_d(\imi t)| \, dt
    \ll L^{c_0} \sum_{d \le N_3} Xd^{-1} + \tau(m)TX^{11/20}L^{c_0} 
  \end{align}
  Combining \eqref{2.1.a}, \eqref{2.1.b}, and \eqref{2.1.d} we obtain
  \[
    \int_{-T}^T |F(\imi t, \chi_0)| \, dt \ll \big( X + m^{0.01}TX^{11/20} \big) L^{c_0 + 1}.
  \]
  
  \medskip
  
  \paragraph{\emph{Step 2}}
  Next, we treat the case where $\max(N_1, N_2) \ge X^{9/20}$. Suppose first that $X^{9/20} \ll N_1 \ll X^{11/20}$. By \cite[Theorem 9.12]{IK} and \eqref{1.2},
  \begin{equation}\label{2.1.1}
    \sum_{\chi \in \mathcal H(m, r, Q)} \int_{-T}^T |F_1(\imi t, \chi)|^2 \, dt 
    \ll (N_1 + H)N_1L^{\kappa^2 + 2}.
  \end{equation}
  Similarly, 
  \begin{equation}\label{2.1.2}
    \sum_{\chi \in \mathcal H(m, r, Q)} \int_{-T}^T |\tilde F_2(\imi t, \chi)|^2 \, dt 
    \ll (N_2N_3 + H)N_2N_3L^{\nu^2 + 2\nu + 3},
  \end{equation}
  where $\tilde F_2(s, \chi) = F_2(s, \chi)F_3(s, \chi)$ is a Dirichlet polynomial with coefficients $\tilde b_2(n)$ subject to
  \[
    \big| \tilde b_2(n) \big| \le \sum_{n = uv} \tau_{\nu}(u) \le \tau_{\nu + 1}(n).
  \]
  Using \eqref{2.1.1}, \eqref{2.1.2}, and the Cauchy--Schwarz inequality, we find that the left side of \eqref{1.3} is bounded above by  
  \begin{align*}
    & \big( X^{1/2} + (HN_1)^{1/2} + (HN_2N_3)^{1/2} + H \big) X^{1/2}L^{c_1}  \\
    \ll \; & \big( X + H^{1/2}X^{31/40} + HX^{1/2} \big) L^{c_1}  
    \ll \big( X + HX^{11/20} \big) L^{c_1},
  \end{align*}
  where $c_1 = c_1(\kappa, \nu) = \kappa^2 + \nu^2 + 4$. Since an obvious modification of this argument establishes \eqref{1.3} when $N_2 \gg X^{9/20}$, we may assume for the remainder of the proof that 
  \begin{equation}\label{3.09}
    \max(N_1, N_2) \le X^{9/20}.
  \end{equation}
  
  \medskip
  
  \paragraph{\emph{Step 3}}
  Suppose that hypothesis (ii) holds. By a standard argument,
  \begin{equation}\label{2.34}
    \sum_{\chi \in \mathcal H(m, r, Q)} \int_{-T}^T \left| F (\imi t, \chi) \right| \, dt
    \ll \sum_{j = 1}^R \left| F (\imi t_j, \chi_j) \right|,
  \end{equation}
  where $(t_1, \chi_1), \dots, (t_R, \chi_R)$ are well-spaced points. The points $(t_j, \chi_j)$ such that
  \[
    F_i(\imi t_j, \chi_j) \ll X^{-1} \qquad \text{for some } i = 1, 2, 3
  \]
  contribute at most
  \[
    RX^{-1}X^{1.01} \ll RX^{0.01} \ll HX^{0.01}
  \]
  to the right side of \eqref{2.34}. We divide the remaining points $(t_j, \chi_j)$ into $O\big( L^3 \big)$ subsets so that for the points in a particular subset $\mathcal S(V_1, V_2, V_3)$ we have
  \begin{equation}\label{3.05}
    V_i \le |F_i(\imi t_j, \chi_j)| \le 2V_i \qquad (i = 1, 2, 3).
  \end{equation}
  We obtain that 
  \begin{equation}\label{3.06}
    \sum_{\chi \in \mathcal H(m, r, Q)} \int_{-T}^T \big| F (\imi t, \chi) \big| \, dt
    \ll L^3V_1V_2V_3|\mathcal S(V_1, V_2, V_3)| + HX^{0.01}
  \end{equation}
  for some $V_1, V_2, V_3$ subject to
  \begin{equation}\label{3.064}
    X^{-1} \le V_i \le N_iL^{\kappa + \nu}.
  \end{equation}
  Thus, it suffices to show that 
  \begin{equation}\label{3.065}
    V_1V_2V_3|\mathcal S(V_1, V_2, V_3)| \ll \big( X + HX^{11/20} \big) L^{c_2 + \kappa + \nu},
  \end{equation} 
  where $c_2 = c_2(\kappa, \nu) = 3\max(\kappa^2, \nu^2) + 15$. To derive this bound, we apply Lemma \ref{l2.2} to $F_1(s, \chi)$, $F_2(s, \chi)$, and $F_3(s, \chi)^2$ and find that
  \begin{align}\label{3.010}
    |\mathcal S(V_1, V_2, V_3)| \ll \min \big\{ 
    & N_1^2V_1^{-2} + HN_1\min \big( V_1^{-2}, N_1^3V_1^{-6} \big), \\
    & N_2^2V_2^{-2} + HN_2\min \big( V_2^{-2}, N_2^3V_2^{-6} \big), \notag\\
    & N_3^4V_3^{-4} + HN_3^2\min \big( V_3^{-4}, N_3^6V_3^{-12} \big) \big\} L^{c_2}. \notag
  \end{align}

  \medskip
  
  \paragraph{\emph{Step 4}}
  Suppose that hypothesis (i) holds and $\mathcal H(m, r, Q)$ contains no principal characters. We combine the argument from Step 3 with the observation that under the present assumptions we also have the estimate
  \[
    |\mathcal S(V_1, V_2, V_3)| \ll HN_3^2V_3^{-4}L^{10}
  \] 
  (this follows from \eqref{3.05} and Lemma \ref{l2.3}). Thus, we obtain \eqref{3.06} with 
  \begin{align}\label{3.011}
    |\mathcal S(V_1, V_2, V_3)| \ll &\min \big\{ 
    N_1^2V_1^{-2} + HN_1 \min \big( V_1^{-2}, N_1^3V_1^{-6} \big), \\
    & \, N_2^2V_2^{-2} + HN_2 \min \big( V_2^{-2}, N_2^3V_2^{-6} \big), HN_3^2V_3^{-4} \big\} L^{c_2}.\notag
  \end{align}
  
  We must also supply a proof of the bound \eqref{2.1.c} used in Step 1. In this case we have to deal with well-spaced points $(t_1, \chi^0), \dots, (t_R, \chi^0)$, where $|t_j| \ge M_d$ and $\chi^0$ is the trivial character: $\chi^0(n) = 1$ for all $n$. Thus, Lemma \ref{l2.3} (with $H = T$) can again be used to show that
  \begin{align}\label{3.011a}
    |\mathcal S(V_1, V_2, V_3)| \ll &\min \big\{ 
    N_1^2V_1^{-2} + TN_1 \min \big( V_1^{-2}, N_1^3V_1^{-6} \big), \\
    & \, N_2^2V_2^{-2} + TN_2 \min \big( V_2^{-2}, N_2^3V_2^{-6} \big), TM_d^2V_3^{-4} \big\} L^{c_2}.\notag
  \end{align}
  
  \medskip
  
  \paragraph{\emph{Step 5}}
  The remainder of the proof is a case-by-case analysis that derives \eqref{3.065} from \eqref{3.09}, \eqref{3.064}, and \eqref{3.011} under hypothesis (i) and from \eqref{3.09}, \eqref{3.064}, and   \eqref{3.010} under hypothesis (ii). We write
  \[
    \Gamma_i = N_i^2V_i^{-2}, \qquad \Delta_i = \min \big( V_i^{-2}, N_i^3V_i^{-6} \big), \qquad \Delta_i(\alpha) = N_i^{3\alpha}V_i^{-2 - 4\alpha}
  \]
  and remark that $\Delta_i \le \Delta_i(\alpha)$ for all $0 \le \alpha \le 1$.
  
  \medskip
  
  \paragraph{\em Case 1:} $\Gamma_1 \ge HN_1\Delta_1$ and $\Gamma_2 \ge HN_2\Delta_2$. Then, by \eqref{3.064} and \eqref{3.010} or \eqref{3.011},
  \begin{align*}
    V_1V_2V_3|\mathcal S(V_1, V_2, V_3)| 
    &\ll V_1V_2V_3 \min \big\{ \Gamma_1, \Gamma_2 \big\} L^{c_2} \\
    &\ll V_1V_2V_3 (\Gamma_1\Gamma_2)^{1/2} L^{c_2} \\
    &\ll N_1N_2V_3L^{c_2} \ll XL^{c_2 + \kappa + \nu}.
  \end{align*}
  
  \smallskip
  
  \paragraph{\em Case 2:} $\Gamma_1 \le HN_1\Delta_1$, $\Gamma_2 \le HN_2\Delta_2$, and $\Gamma_3^2 \ge HN_3^2\Delta_3^2$. This case occurs only when we argue under hypothesis (ii). By \eqref{3.010} and the hypothesis $N_3 \le X^{8/35}$, we get
  \begin{align*}
    V_1V_2V_3|\mathcal S(V_1, V_2, V_3)| &\ll 
    V_1V_2V_3 \min \big\{ HN_1\Delta_1, HN_2\Delta_2, \Gamma_3^2 \big\} L^{c_2} \\
    &\ll V_1V_2V_3 \big( HN_1\Delta_1(1/6) \big)^{3/8} \big( HN_2\Delta_2(1/6) \big)^{3/8} \Gamma_3^{1/2} L^{c_2}\\
    &\ll H^{3/4}(X^9N_3^7)^{1/16} L^{c_2} \ll \big( X + HX^{11/20} \big) L^{c_2},    
  \end{align*}
  where the last step uses that 
  \[
    H^{3/4}(X^9N_3^7)^{1/16} \ll H^{3/4} X^{53/80}
    = X^{1/4} \big(HX^{11/20}\big)^{3/4}.
  \] 
   
  \smallskip
  
  \paragraph{\em Case 3:} $\Gamma_1 \le HN_1\Delta_1$, $\Gamma_2 \le HN_2\Delta_2$, and $\Gamma_3^2 \le HN_3^2\Delta_3^2$. When $N_3 \le X^{1/5}$, \eqref{3.010} yields
  \begin{align*}
    V_1V_2V_3|\mathcal S(V_1, V_2, V_3)| &\ll 
    V_1V_2V_3 \min \big\{ HN_1\Delta_1, HN_2\Delta_2, HN_3^2\Delta_3^2 \big\} L^{c_2} \\
    &\ll HV_1V_2V_3 (N_1\Delta_1(1/22)N_2\Delta_2(1/22))^{11/24}(N_3\Delta_3(1))^{1/6} L^{c_2} \\ 
    &\ll H(X^{25}N_3^7)^{1/48} L^{c_2} \ll HX^{11/20}L^{c_2}.
  \end{align*}
  On the other hand, when $N_3 \ge X^{1/5}$, both \eqref{3.010} and \eqref{3.011} yield
  \begin{align*}
    V_1V_2V_3 |\mathcal S(V_1, V_2, V_3)| 
    &\ll HV_1V_2V_3 \big( N_1\Delta_1(1/6) N_2\Delta_2(1/6) \big)^{3/8} (N_3\Delta_3(0))^{1/2} L^{c_2}\\
    &\ll H(X^9N_3^{-1})^{1/16} L^{c_2} \ll HX^{11/20} L^{c_2}.    
  \end{align*}
    
  \smallskip
  
  \paragraph{\em Case 4:} $\Gamma_1 \ge HN_1\Delta_1$, $\Gamma_2 \le HN_2\Delta_2$, and $\Gamma_3^2 \ge HN_3^2\Delta_3^2$. Again, this only occurs when we argue from \eqref{3.010}. By \eqref{3.09}, \eqref{3.010}, and the hypothesis $N_3 \le X^{8/35}$, 
  \begin{align*}
    V_1V_2V_3 |\mathcal S(V_1, V_2, V_3)| 
    &\ll V_1V_2V_3 \min \big\{ \Gamma_1, HN_2\Delta_2, \Gamma_3^2 \big\} L^{c_2} \\
    &\ll V_1V_2V_3 (\Gamma_1\Gamma_3)^{1/2}(HN_2\Delta_2(1/2))^{1/4} L^{c_2} \\ 
    &\ll H^{1/4}X^{5/8}(N_1N_3)^{3/8} L^{c_2} \ll \big( X + HX^{11/20} \big) L^{c_2},
  \end{align*}
  where the last step uses that 
  \[
    H^{1/4}X^{5/8}(N_1N_3)^{3/8} \ll H^{1/4}X^{197/224} \ll X^{3/4} \big(HX^{11/20}\big)^{1/4}.
  \]

  \smallskip
  
  \paragraph{\em Case 5:} $\Gamma_1 \ge HN_1\Delta_1$, $\Gamma_2 \le HN_2\Delta_2$, and $\Gamma_3^2 \le HN_3^2\Delta_3^2$. When $N_3 \le X^{1/5}$, \eqref{3.010} yields
  \begin{align*}
    V_1V_2V_3 |\mathcal S(V_1, V_2, V_3)|  
    &\ll V_1V_2V_3 \min \big\{ \Gamma_1, HN_2\Delta_2, HN_3^2\Delta_3^2 \big\} L^{c_2} \\
    &\ll V_1V_2V_3 \Gamma_1^{1/2}(HN_2\Delta_2(1/10))^{5/12}(HN_3^2\Delta_3(1)^2)^{1/12}L^{c_2} \\ 
    &\ll H^{1/2}(N_1^{11}N_3^3)^{1/24}X^{13/24} L^{c_2} 
    \ll \big( X + HX^{11/20} \big) L^{c_2},
  \end{align*}
  where the last step uses that 
  \[
    H^{1/2}(N_1^{11}N_3^3)^{1/24}X^{13/24} \ll H^{1/2} X^{371/480} 
    \ll X^{1/2} \big( HX^{11/20} \big)^{1/2}.
  \]
  On the other hand, when $N_3 \ge X^{1/5}$, by \eqref{3.010} or \eqref{3.011}, 
  \begin{align*}
    V_1V_2V_3 |\mathcal S(V_1, V_2, V_3)| 
    &\ll V_1V_2V_3 \Gamma_1^{1/2}(HN_2\Delta_2(1/2))^{1/4}(HN_3^2\Delta_3(0)^2)^{1/4}L^{c_2} \\ 
    &\ll H^{1/2}(X^5N_1^3N_3^{-1})^{1/8} L^{c_2} 
    \ll \big( X + HX^{11/20} \big) L^{c_2},
  \end{align*}
  because
  \[
    H^{1/2}(X^5N_1^3N_3^{-1})^{1/8} \ll H^{1/2} X^{123/160} \ll X^{1/2} \big( HX^{11/20} \big)^{1/2}.
  \]
 
  \smallskip
  
  \paragraph{\em Case 6:} $\Gamma_1 \le HN_1\Delta_1$ and $\Gamma_2 \ge HN_2\Delta_2$. This case can be split into two subcases that can be handled similarly to Cases 4 and 5.
\end{proof}

We conclude this section with a technical lemma, which will be needed in the next section.

\begin{lemma}\label{l3.1}
  Suppose that $2 \le T \le M < N$ and $f: \mathbb N^2 \to \mathbb C$ is a function such that the inequality
  \begin{equation}\label{3.1.1}
    \sum_m \int_{-U}^U \bigg| \sum_{n \le N} f(m, n) n^{\imi t} \bigg| \, dt \le A + BU
  \end{equation}
  holds for all $U \ge 2$. Then
  \begin{equation}\label{3.1.2}
    \sum_m \int_{-T}^T \bigg| \sum_{n \le M} f(m, n) n^{\imi t} \bigg| \, dt \ll (A + BT)\log^2 N.
  \end{equation}
\end{lemma}

\begin{proof}
  Let $g$ denote the indicator function of $[-T, T]$ and let $h$ be the function constructed in \cite[Lemma 13.11]{IK} with $z = N$. Then
  \begin{equation}\label{3.1.3}
    |h(u)| \ll \min \big\{ \log N, |u|^{-1}, N|u|^{-2} \big\}
  \end{equation}
  and 
  \[
    \int_{-\infty}^{\infty} h(u) \left( \frac mn \right)^{\imi u} du = \begin{cases}
      1 & \text{if } m \le n, \\
      0 & \text{if } m > n,
    \end{cases}
  \]
  for any pair of integers $m, n$ such that $1 \le m, n \le N$. Thus, 
  \[
    \sum_{n \le M} f(m, n) n^{\imi t} = \int_{-\infty}^{\infty} 
    \bigg\{ \sum_{n \le N} f(m, n) n^{\imi (t + u)} \bigg\} h(u)M^{-\imi u} \, du,
  \]
  assuming (as we may) that $M$ is an integer. It follows that the left side of \eqref{3.1.2} does not exceed
  \begin{align*}
    & \; \sum_m \int_{-\infty}^{\infty} g(t) \int_{-\infty}^{\infty} |h(u)| 
    \bigg| \sum_{n \le N} f(m, n) n^{\imi (t + u)} \bigg| \, dudt \\
    =& \; \sum_m \int_{-\infty}^{\infty} \bigg| \sum_{n \le N} f(m, n) n^{\imi \tau} \bigg| 
    \bigg\{ \int_{-\infty}^{\infty} g(\tau - u)|h(u)| \, du \bigg\} d\tau \\
    \ll& \; T(\log N) \sum_m \int_{-\infty}^{\infty} \bigg| \sum_{n \le N} f(m, n) n^{\imi \tau} \bigg| 
    \min \big\{ T^{-1}, |\tau|^{-1}, N|\tau|^{-2} \big\} d\tau,
  \end{align*}
  where the last step uses \eqref{3.1.3} and the definition of $g$. The desired conclusion now follows by a standard dyadic argument.
\end{proof}

\section{Proof of Theorem \ref{th1}}

In this section we deduce Theorem \ref{th1} from Theorem \ref{th2.1} and Heath-Brown's identity for $\Lambda(n)$. We apply Heath-Brown's identity in the following form (see \cite[Lemma 1]{HB} or \cite[Proposition 13.3]{IK} with $k = 10$): if $n \le x$, then
\begin{equation}\label{4.1}
  \Lambda(n) = \sum_{j = 1}^{10} \binom {10}j (-1)^j
  \sum_{\substack{ n = m_1 \cdots m_{2j}\\
  m_1, \dots, m_j \le x^{1/10}}}
  \mu(m_1) \cdots \mu(m_j) \log m_{2j}.
\end{equation}
By \eqref{4.1} with $x = 2N$ and a standard splitting argument,
\[
  \sum_{N < n \le 2N} \Lambda(n) \chi(n) n^{-s} \ll
  \sum_{\mathbf M} \bigg| \sum_{N < n \le 2N} a(n; \mathbf M) \chi(n) n^{-s} \bigg|,
\]
where $\mathbf M$ runs over $O( L^{19} )$ vectors $\mathbf M = (M_1, \dots, M_{2j})$, $j \le 10$, subject to
\[
  M_1, \dots, M_j \ll N^{1/10}, \qquad N \ll M_1 \cdots M_{2j} \ll N,
\]
and
\[
  a(n; \mathbf M) = \sum_{\substack{ n = m_1 \cdots m_{2j}\\ M_i < m_i \le 2M_i}} \mu(m_1) \cdots \mu(m_j)(\log m_{2j}).
\]
Thus, the left side of \eqref{i.2} is bounded above by
\[
  L^{19} \sum_{\chi \in \mathcal H(m, r, Q)} \int_{-T}^T \bigg| \sum_{N < n \le 2N} a(n; \mathbf M) \chi(n) n^{-\imi t} \bigg| \, dt
\]
for some fixed choice of $\mathbf M$ as above. Thus, if we show that
\begin{equation}\label{4.2}
  \sum_{\chi \in \mathcal H(m, r, Q)} \int_{-T}^T \bigg| \sum_n a(n; \mathbf M) \chi(n) n^{-\imi t} \bigg| \, dt
  \ll \big( N + HN^{11/20} \big) L^{1020},
\end{equation}
the desired result (with $C = 1100$) will follow by Lemma \ref{l3.1}.

The Dirichlet polynomial on the right side of \eqref{4.2} is the product of $2j$, $j \le 10$, Dirichlet polynomials of the form \eqref{i.1} with coefficients $a_n = \mu(n)$, $a_n = 1$, or $a_n = \log n$. Furthermore, the single logarithmic weight can be removed by partial summation. Therefore, we may assume that 
\[
  a(n; \mathbf M) = L \sum_{\substack{ n = m_1 \cdots m_{2j}\\ M_i < m_i \le M_i'}} \mu(m_1) \cdots \mu(m_j),
\]
where $M_i < M_i' \le 2M_i$ (in reality, $M_i' = 2M_i$ except for $i = 2j$). We may now assume that $M_{j + 1} \le \cdots \le M_{2j}$. We proceed to show that 
\[
  a(n; \mathbf M) = L \sum_{n = n_1n_2n_3} b_1(n_1)b_2(n_2)b_3(n_3),
\]
where the coefficients on the right yield a Dirichlet polynomial \eqref{1.1} that satisfies at least one of the hypotheses (i) or (ii) of Theorem \ref{th2.1}. The analysis involves several cases depending on the sizes of $M_1, \dots, M_{2j}$.

\medskip
  
\paragraph{\em Case 1:} $M_{2j} \gg N^{9/20}$. Assuming that $j \ge 2$ (the case $j = 1$ is similar and easier), we group the variables $m_1, \dots, m_{2j}$ into $n_1, n_2, n_3$ as follows:
\[
  n_1 = m_3 \cdots m_{2j - 1}, \qquad n_2 = m_1m_2, \qquad n_3 = m_{2j}.
\]
Since $M_1 \cdots M_{2j - 1} \ll N^{11/20}$, this yields a polynomial $F(s, \chi)$ satisfying hypothesis (i) of Theorem \ref{th2.1}.

\smallskip
  
\paragraph{\em Case 2:} $M_{2j} \ll N^{9/20} \ll M_1 \cdots M_jM_{2j}$. Let $i$ be the least integer for which $M_1 \cdots M_iM_{2j} \gg N^{9/20}$. Since $M_i \ll N^{1/10}$, we have
\[
    N^{9/20} \ll M_1 \cdots M_iM_{2j} \ll N^{11/20}.
\]
Hence, the choice
\[
  n_1 = m_1 \cdots m_im_{2j}, \qquad n_2 = m_{i + 1} \cdots m_{2j - 1}, \qquad n_3 = 1
\]
yields an $F(s, \chi)$ that satisfies hypothesis (ii) of Theorem \ref{th2.1}.
  
\smallskip
  
\paragraph{\em Case 3:} $M_1 \cdots M_jM_{2j} \ll N^{9/20}$. Let $\ell$ be the least positive integer such that
\[
    M_1 \cdots M_j M_{\ell} \cdots M_{2j} \ll N^{9/20}.
\]
We consider three subcases.

\smallskip
  
\paragraph{\em Case 3.1:} $M_{\ell - 1} \cdots M_{2j} \ll N^{11/20}$. Then we can argue similarly to Case 2 to find an $i$, $0 \le i \le j$, for which
\[
    N^{9/20} \ll M_1 \cdots M_i M_{\ell - 1} \cdots M_{2j} \ll N^{11/20}.
\]
Again, we will have $F(s, \chi)$ that satisfies hypothesis (ii) of Theorem \ref{th2.1}.

\smallskip
  
\paragraph{\em Case 3.2:} $M_{\ell - 1} \cdots M_{2j} \gg N^{11/20}$ and $M_{\ell - 1} \ll N^{8/35}$. Then we define
\[
  n_1 = m_1 \cdots m_jm_{\ell} \cdots m_{2j}, \qquad n_2 = m_{j + 1} \cdots m_{\ell - 2}, \qquad  n_3 = m_{\ell - 1}.
\]
Since $M_{j + 1} \cdots M_{\ell - 2} \ll N^{9/20}$, we again get an $F(s, \chi)$ that satisfies hypothesis (ii) of Theorem \ref{th2.1}.
  
\smallskip
  
\paragraph{\em Case 3.3:} $M_{\ell - 1} \cdots M_{2j} \gg N^{11/20}$ and $M_{\ell - 1} \gg N^{8/35}$. This may occur only with $\ell = 2j$. Then 
\[
    M_1 \cdots M_{2j - 2} \ll NM_{2j - 1}^{-2} \ll N^{19/35} \ll N^{11/20} 
    \quad \text{and} \quad M_{2j - 1} \ll M_{2j} \ll N^{9/20}.
\]
We write
\[
  b_1(n) = \sum_{n = m_1 \cdots m_{2j - 2}} \mu(m_1) \cdots \mu(m_j), \qquad n_2 = m_{2j - 1}, \qquad n_3 = m_{2j},
\]
and we obtain an $F(s, \chi)$ that satisfies hypothesis (i) of Theorem \ref{th2.1}.

\medskip

The desired bound \eqref{4.2} follows on noting that the arising coefficients satisfy \eqref{1.2} with $\kappa, \nu$ for which $c(\kappa, \nu) \le c(18, 2) = 1012$.

\section{Exponential sums twisted by characters}
\label{s5}

In this section we estimate the exponential sum
\begin{equation}\label{1.6}
  W(\beta, \chi) = \sum_{N < p \le 2N} (\log p) \chi(p) e \big( \beta p^k \big),
\end{equation}
where $k$ is a positive integer, $\beta$ is ``small'', and $\chi$ is Dirichlet character. Such exponential sums arise in dealings with the major arcs in the Waring--Goldbach problem and related questions. In particular, in the proof of Theorem \ref{th3} we need the case $k = 1$ of our estimates.

\begin{lemma}\label{lemma 1.3}
  Suppose that $N \ge 2$ and $0 \le \Delta \le N^{1 - k}$. Suppose also that $\mathcal H(m, r, Q)$ is a set of characters as in Theorem \ref{th1} and $W(\beta, \chi)$ is defined by \eqref{1.6}. Then
  \begin{equation}\label{1.8}
    \sum_{\chi \in \mathcal H(m, r, Q)} \max_{\Delta \le |\beta| \le 2\Delta} |W(\beta, \chi)|
    \ll T_0^{-1/2}L^{C + 1}\big( N + HN^{11/20} \big),
  \end{equation}
  where $T_0 = 1 + \Delta N^k$, $H = mr^{-1}Q^2T_0$, $L = \log N$, and $C$ is the constant appearing in Theorem \ref{th1}.
\end{lemma}

\begin{proof}
We first replace $W(\beta, \chi)$ by the exponential sum
\[
  \tilde W(\beta, \chi) = \sum_{N < n \le 2N} \Lambda(n)\chi(n)e \left( \beta n^k \right)
\]
using that
\begin{equation}\label{1.9}
W(\beta, \chi) = \tilde W(\beta, \chi) + O \big( N^{1/2} \big).
\end{equation}
By Perron's formula \cite[Proposition 5.54]{IK}, for $N < M \le 2N$,
\begin{equation}\label{1.14}
  \sum_{N < n \le M} \Lambda(n)\chi(n) = \frac 1{2\pi \imi} \int_{b - \imi T_1}^{b + \imi T_1} F(s, \chi) \frac {M^s - N^s}s \, ds + O \left( \frac {NL^2}{1 + T_1\| M \|} \right),
\end{equation}
where $0 < b < (\log N)^{-1}$, $T_1 = (HN)^{10}$, $\| M \|$ is the distance from $M$ to the nearest integer, and
\[
  F(s, \chi) = \sum_{N < n \le 2N} \Lambda(n) \chi(n) n^{-s}.
\]
Hence, by partial summation, 
\begin{equation}\label{1.15}
  \tilde W(\beta, \chi) = \frac 1{2\pi \imi } \int_{b - \imi T_1}^{b + \imi T_1} F(s, \chi) V(s, \beta) \, ds + O(1),
\end{equation}
where 
\[
  V(s, \beta) = \int_N^{2N} y^{s - 1} e \big( \beta y^k \big) dy.
\]
By \cite[Lemma 8.10]{IK}, for $\Delta \le |\beta| \le 2\Delta$,
\begin{equation}\label{1.16}
  V(\sigma + \imi t, \beta) \ll N^{\sigma} \min \big\{ T_0^{-1/2}, \sup_{N \le y \le 2N} |t + 2k\pi\beta y^k|^{-1} \big\},
\end{equation}
Combining \eqref{1.15} and \eqref{1.16} and letting $b \downarrow 0$, we obtain
\[
  \tilde W(\beta, \chi) \ll T_0^{1/2} \int_{-T_1}^{T_1} \left| F (\imi t, \chi) \right|
  \frac {dt}{T_0 + |t|} + 1.
\]
Recalling \eqref{1.9}, we deduce that the right side of \eqref{1.8} is bounded above by
\begin{equation}\label{1.18}
  LT_0^{1/2}T^{-1} \sum_{ \chi \in \mathcal H(m, r, Q)} \int_{-T}^T \big| F(\imi t, \chi) \big| \, dt 
  + |\mathcal H|N^{1/2},
\end{equation}
for some $T$ in the range $T_0 \le T \le T_1$. The desired result now follows from \eqref{i.2}.
\end{proof}

We now define the exponential integral
\begin{equation}\label{1.7}
  v(\beta; X) = \int_X^{2X} e \big( \beta y^k \big) \, dy .
\end{equation}

\begin{lemma}\label{lemma 1.4}
Suppose that $N \ge 2$, $1 \le Q \le N$, and $0 \le \Delta \le N^{1 - k}$. 
Let $W(\beta; \chi)$ be defined by \eqref{1.6}. Then, for any fixed $A > 0$ and $\delta > 0$,
\begin{equation}\label{1.20}
  \sum_{Q < q \le 2Q} \; \sideset{}{^*} \sum_{\chi\!\!\!\!\! \mod q} 
  \max_{\Delta \le |\beta| \le 2\Delta}  |W(\beta, \chi)| \ll 
  NQ^{\delta}L^{-A} + Q^2T_0^{1/2}N^{11/20}L^{C + 1},
\end{equation}
where $T_0 = 1 + \Delta N^k$, $L = \log N$, and $C$ is the constant appearing in Theorem \ref{th1}.
Furthermore, for any fixed $A > 0$ we have
\begin{equation}\label{1.21}
  W(\beta, \chi^0) - v(\beta; N) \ll NL^{-A} + T_0^{1/2}N^{11/20}L^{C + 1},
\end{equation}
where $v(\beta; N)$ is defined by \eqref{1.7} and $\chi^0$ is the trivial character. 
In both \eqref{1.20} and \eqref{1.21} the implied constants may depend on $A$, 
and the implied constant in \eqref{1.20} may also depend on $\delta$.
\end{lemma}

\begin{proof}
The first claim follows from Lemma \ref{lemma 1.3} and the Siegel--Walfisz theorem 
in the form of \cite[(5.79)]{IK}. Put $B = (2 + \delta^{-1})(A + C + 1)$. 
If $Q \ge L^B$ or $\Delta \ge L^BN^{-k}$, we have
\[
  NT_0^{-1/2}L^{C + 1} \ll XQ^{\delta}L^{-A}
\]
and \eqref{1.20} follows from \eqref{1.8} with $m = r = 1$. On the other hand, 
if $Q \le L^B$ and $\Delta \le L^BN^{-k}$, we find by partial summation that 
the left side of \eqref{1.20} is bounded above by
\[
  L^{3B + 1} \max_{N < N_1 \le 2N} \bigg| \sum_{N < p \le N_1} \chi(p) \bigg| \ll NL^{-A},
\]
by the aforementioned version of the Siegel--Walfisz theorem.

The proof of the second claim is similar, except that it appeals to the case 
$m = r = Q = 1$ of Lemma  \ref{lemma 1.3} and to the prime number theorem 
(which is why we need to include the term $v(\beta; N)$ on the left side of 
\eqref{1.21}).
\end{proof}

\begin{lemma}
\label{lemma 1.5}
  Suppose that $N \ge 2$ and $N^{-k} \le \Delta \le N^{1 - k}$. Suppose also that $\mathcal H(m, r, Q)$ is a set of characters as in Theorem \ref{th1} and $W(\beta, \chi)$ is defined by \eqref{1.6}. Then
  \begin{equation}\label{1.22}
    \sum_{\chi \in \mathcal H(m, r, Q)} \bigg\{ \int_{-\Delta}^{\Delta} |W(\beta, \chi)|^2 \, d\beta \bigg\}^{1/2}
    \ll N^{-k/2}L^{C + 1} \big( N + HN^{11/20} \big),
  \end{equation}
  where $H = mr^{-1}Q^2\Delta N^k$, $L = \log N$, and $C$ is the constant from Theorem \ref{th1}.
\end{lemma}

\begin{proof}
By \cite[Lemma 1.9]{Mo}, we have
\begin{align}\label{1.23}
  \int_{-\Delta}^{\Delta} |W(\beta, \chi)|^2 \, d\beta &\ll \Delta^2 
  \int_{-\infty}^{\infty}\bigg| \sum_{u(y) < p \le v(y)} (\log p)\chi(p) \bigg|^2 dy \\
  &\ll \Delta^2X^k \bigg| \sum_{M < n \le M + Y} \Lambda(n) \chi(n) \bigg|^2 + \Delta^2X^{k + 1}, \notag
\end{align}
where $u(y) = \max (N, y^{1/k})$, $v(y) = \min(2N, (y + (2\Delta)^{-1})^{1/k})$, and
\begin{equation}\label{1.24}
  N < M \le 2N, \qquad Y \ll \Delta^{-1} N^{1 - k}.
\end{equation}
Without loss of generality, we may assume that the distance from $M$ to the 
nearest integer is $1/2$ and that $Y$ is an integer. We then appeal to Perron's formula to derive
\[
\sum_{M < n \le M + Y} \Lambda(n) \chi(n)
\ll \bigg| \int_{b - \imi T_1}^{b + \imi T_1}
F(s, \chi) \frac {(M + Y)^s - M^s}s \, ds \bigg| + 1,
\]
where $0 < b < L^{-1}$, $T_1 = (HN)^{10}$, and $F(s, \chi)$ is the Dirichlet polynomial appearing in the proof of Lemma \ref{lemma 1.3}. Hence, as in that proof, 
\begin{equation}\label{1.25}
  \sum_{M < n \le M + Y} \Lambda(n) \chi(n)
  \ll \int_{-T_1}^{T_1} \big| F(\imi t, \chi) \big| \, \frac {dt}{T_0 + |t|} + 1,
\end{equation}
where $T_0 = \Delta N^k$. By \eqref{1.23} and \eqref{1.25}, the left side of \eqref{1.22} is bounded above by
\[
  \Delta N^{k/2}LT^{-1} \sum_{\chi \in \mathcal H(m, r, Q)} \int_{-T}^T \big| F(\imi t, \chi) \big| \, dt + HN^{(1 - k)/2},
\]
where $T$ is subject to $T_0 \le T \le T_1$. The desired result now follows from \eqref{i.2}.
\end{proof}

\section{Proof of Theorem \ref{th3}}

Since the proof follows closely the proof of the main result in \cite{jLT}, we only describe the necessary changes. Let $N$ be a large parameter chosen as in \cite[Lemma 2.3]{jLT} and set
\begin{equation}\label{6.1}
  P = (N/B)^{9/20}, \qquad L = \log N, \qquad Q = N / (PL^2).
\end{equation}
We note that the improvement on the result of Liu and Tsang arises from the choice of $P$ in \eqref{6.1}: the respective choice in \cite{jLT} is $P = (N/B)^{2/5}$ (see \cite[(2.1)]{jLT}). In order to justify the analysis in \cite{jLT} for this larger value of $P$, we must establish appropriate variants of \cite[Lemmas 3.2 and 3.3]{jLT}. 

Let $N_j = N/|a_j|$, $N^{1/10} \le R \le P$, and $g, D$ be positive integers. Define
\[
  K_j(g; R) = \sum_{R < r \le 2R} \frac{\sqrt{([g, r], D)}}{[g, r]} \sideset{}{^*}\sum_{\chi \modulo r}
  \bigg( \int_{-1/(RQ)}^{1/(RQ)} \big| W_j( a_j\beta; \chi ) \big|^2 \, d\beta \bigg)^{1/2},
\]
where $W_j(\beta; \chi)$ is the sum \eqref{1.6} with $N = N_j$ and $k = 1$. In order to prove \cite[Lemma 3.2]{jLT} with $P$ as in \eqref{6.1}, we need to show that
\begin{equation}\label{6.2}
  K_j(g; R) \ll g^{-1} \sqrt{(g, D)}\tau(gD)^2 N_jN^{-1/2}L^c
\end{equation}
for some absolute constant $c$. By \cite[(5.20)]{jLT}, 
\begin{equation}\label{6.3}
  K_j(g; R) \ll \frac {\sqrt{(g, D)}}{gR} 
  \sum_{\substack{d \mid gD\\ d \le 2R}} d\tau(d) \tilde K_j(d; R),
\end{equation}
where 
\[
  \tilde K_j(d; R) = \sum_{\chi \in \mathcal H(1, d, 2R)} 
  \bigg( \int_{-1/(RQ)}^{1/(RQ)} \big| W_j( a_j\beta; \chi ) \big|^2 \, d\beta \bigg)^{1/2}.
\]
By Lemma \ref{lemma 1.5} with $k = 1$, 
\begin{align*}
  \tilde K_j(d; R) &\ll |a_j|^{-1/2} \sum_{\chi \in \mathcal H(1, d, 2R)} 
  \bigg( \int_{-|a_j|/(RQ)}^{|a_j|/(RQ)} \big| W_j( \beta; \chi ) \big|^2 \, d\beta \bigg)^{1/2} \\
  &\ll N^{-1/2}L^{C + 1} \big( N_j + H_jN_j^{11/20} \big),
\end{align*}
where $C$ is the constant appearing in Theorem \ref{th1} and 
\[
  H_j = d^{-1}R^2 \big( |a_j|/ (RQ) \big) N_j \ll d^{-1}PRL^2 \ll d^{-1}RN_j^{9/20}L^2.
\]
Thus,
\[
  \tilde K_j(d; R) \ll N_jN^{-1/2}L^{C + 3} \big( R/d + 1 \big).
\]
Clearly, this inequality and \eqref{6.3} imply \eqref{6.2}.

Similarly, we can use Lemma \ref{lemma 1.3} to establish the desired variant of \cite[Lemma 3.3]{jLT}. This completes the proof of the theorem.

\end{document}